\input amstex
\input amsppt.sty
\magnification=\magstep1
\vsize=22.2truecm
\baselineskip=16truept
\nologo
\pageno=1
\topmatter
\def\Z{\Bbb Z}
\def\N{\Bbb N}

\def\bg{\bigg}
\def\({\bg(}
\def\[{\bg[}
\def\){\bg)}
\def\]{\bg]}

\def\mo{\roman{mod}}

\def\bi{\binom}
\def\Str#1#2{\left[\matrix#1\\#2\endmatrix\right]}

\def\eq{\equiv}

\def\bi{\binom}

\def\Proof{\noindent{\it Proof}}

\def\Ack{\medskip\noindent {\bf Acknowledgment}}
\def\pmod #1{\ (\roman{mod}\ #1)}
\def\qbinom #1#2#3{\left[\matrix#1\\#2\endmatrix\right]_{#3}}
\def\fq{\roman{Q}}
\def\eq{\roman{EQ}}
\def\ab #1{\langle{#1}\rangle}
\def\jacob #1#2{\left(\frac{#1}{#2}\right)}
\def\floor #1{\left\lfloor{#1}\right\rfloor}
\def\ceil #1{\left\lceil{#1}\right\rceil}
\title
A congruence involving products of $q$-binomial coefficients
\endtitle
\author
Hao Pan$^1$ and Hui-Qin Cao$^2$
\endauthor
\address
$^1$ Department of Mathematics, Nanjing University,
Nanjing 210093, People's Republic of China
\endaddress
\email{haopan79\@yahoo.com.cn}\endemail
\address
$^2$ Department of Applied Mathematics,
Nanjing Audit University,
Nanjing 210029, People's Public of China
\endaddress
\email{caohq\@nau.edu.cn}\endemail
\abstract In this paper we establish a $q$-analogue of a congruence of Sun concerning the products of binomial coefficients modulo the square of a prime.
\endabstract
\endtopmatter
\document
\TagsOnRight
\heading
1. Introduction
\endheading
In [G], Granville proved the following interesting congruence:
$$
(-1)^{(p-1)(m-1)/2}\prod_{k=1}^{m-1}\bi{p-1}{\left\lfloor kp/m\right\rfloor}\equiv m^p-m+1\ (\mo\ p^2)\tag 1.1
$$
for any prime $p\geq 5$ and $m\geq 2$, where $\floor{x}$ denotes the greatest integer not exceeding $x$.
Later Sun [S] extended Granville's result and showed that
$$
\align
&(-1)^{\frac{p-1}{2}\floor{\frac{m}{2}}}\prod_{1\leq k\leq\floor{m/2}}\binom{p-1}{\floor{pk/m}}\\
\equiv&\cases
\jacob{m}{p}+{\roman{eq}}_p(m)mp\pmod{p^2}\qquad&\text{if }2\nmid m,\\
\jacob{2m}{p}+\jacob{2}{p}{\roman{eq}}_p(m)mp+2\jacob{m}{p}{\roman{eq}}_p(2)p\pmod{p^2}\qquad&\text{if }2\mid m,
\endcases\tag 1.2
\endalign
$$
where $\jacob{\cdot}{p}$ is the Legendre symbol and
$$
{\roman{eq}}_p(m)=\frac{m^{(p-1)/2}-\jacob{m}{p}}{p}
$$
is the Euler quotient.

For an integer $m$ prime to $p$, define the Fermat quotient ${\roman{q}}_p(m)$ by
$$
{\roman{q}}_p(m)=\frac{m^{p-1}-1}{p}.
$$
Observe that
$$
{\roman{q}}_p(2)=\frac{2^{p-1}-1}{p}=\frac{(2^{(p-1)/2}-\jacob{2}{p})(2^{(p-1)/2}+\jacob{2}{p})}{p}\equiv2\jacob{2}{p}{\roman{eq}}_p(2)\pmod{p}.
$$
Then (1.2) can be rewritten as
$$
\align
&(-1)^{\frac{p-1}{2}\floor{\frac{m}{2}}}\jacob{m}{p}\jacob{2}{p}^{m-1}\prod_{1\leq k\leq\floor{m/2}}\binom{p-1}{\floor{pk/m}}\\
\equiv&1+\jacob{m}{p}{\roman{eq}}_p(m)mp+(2\floor{m/2}+1-m){\roman{q}}_p(2)p\pmod{p^2}.\tag 1.3
\endalign
$$

For a non-negative integer $n$, let
$$
[n]_q=\frac{1-q^n}{1-q}=1+q+\cdots+q^{n-1}
$$
and
$$
(x;q)_n=\cases (1-x)(1-xq)\cdots(1-xq^{n-1})\qquad&\text{if }n\geq 1,\\1\qquad&\text{if }n=0.\endcases
$$
And the $q$-binomial coefficients are given by
$$
\qbinom{n}{k}{q}=\frac{(1-q^n)(1-q^{n-1})\cdots(1-q^{n-k+1})}{(1-q^k)(1-q^{k-1})\cdots(1-q)}
$$
for any $k,n\in\N$.
The arithmetic properties of $q$-binomial coefficients have been investigated by serveral authors (e.g., see [A], [C] and [F]).
Recently Pan [P] established a $q$-analogue of
Granville's congruence (1.1). If $p\geq 5$ is a prime and $m\geq 2$ is
an integer with $p\nmid m$, then we have
$$
\align
&(-1)^{(p-1)(m-1)/2}q^{m\sum_{k=1}^{m-1}\bi{\lfloor kp/m\rfloor+1}{2}}\prod_{k=1}^{m-1}\Str{p-1}{\lfloor kp/m\rfloor}_{q^m}\\
\equiv&\frac{m(q^m;q^m)_{p-1}}{(q;q)_{p-1}}-m+1\ (\mo\ [p]_q^2).\tag 1.4
\endalign
$$

In this paper we will give a $q$-analogue of Sun's congruence (1.3).
Suppose that $p$ is an odd prime and $m\geq 2$ is an integer prime to $p$.
It is not difficult to prove that
$$
\frac{(q^m;q^m)_{p-1}}{(q;q)_{p-1}}=\sum_{j=1}^{p-1}\frac{1-q^{jm}}{1-q^j}\equiv1\pmod{[p]_q}.
$$
So we can define the $q$-Fermat quotient by
$$
\fq_p(m,q)=\frac{(q^m;q^m)_{p-1}/(q;q)_{p-1}-1}{[p]_q}.
$$
For any integer $x$, we denote by $\ab{x}_p$ the least non-negative residue of $x$ modulo $p$. Let
$$
R_p(m)=\{1\leq j<p/2:\,\ab{jm}_p>p/2\}.
$$
Then well-known Gauss' lemma asserts that
$$
\jacob{m}{p}=(-1)^{|R_p(m)|}.
$$
Soon we will show that
$$
q^{\sum_{j\in R_p(m)}(p-\ab{jm}_p)}\frac{(q^m;q^m)_{\frac{p-1}{2}}}{(q;q)_{\frac{p-1}{2}}}\equiv\jacob{m}{p}\pmod{[p]_q}.
$$
Define the $q$-Euler quotient by
$$
\eq_p(m,q)=\frac{q^{\sum_{j\in R_p(m)}(p-\ab{jm}_p)}(q^m;q^m)_{\frac{p-1}{2}}\big/(q;q)_{\frac{p-1}{2}}-\jacob{m}{p}}{[p]_{q}}.
$$
\proclaim{Theorem 1.1} Let $p\geq 5$ be a prime and $m\geq 2$ be an integer with $p\nmid m$. Then
$$
\align
&(-1)^{\frac{p-1}{2}\floor{\frac{m}{2}}}\jacob{m}{p}\jacob{2}{p}^{m-1}q^{2m\sum_{k=1}^{\floor{m/2}}\binom{\floor{kp/m}+1}{2}}\prod_{k=1}^{\floor{m/2}}\qbinom{p-1}{\floor{kp/m}}{q^{2m}}\\
\equiv&1+m[p]_q\eq_p^*(m,q)+(2\floor{m/2}+1)[p]_{q^m}\fq_p(2,q^m)-m[p]_q\fq_p(2,q)\\
&+m\(|R_p(m)|+2\sum_{j=1}^{(p-1)/2}\floor{\frac{jm}{p}}\)(1-q^p)\pmod{[p]_q^2},\tag 1.5
\endalign
$$
where
$$
\align
\eq_p^*(m,q)=&\jacob{m}{p}\frac{(1+q^p)\eq_p(m,q^2)}{1+q}\\
=&\frac{\jacob{m}{p}q^{2\sum_{j\in R_p(m)}(p-\ab{jm}_p)}(q^{2m};q^{2m})_{\frac{p-1}{2}}\big/(q^2;q^2)_{\frac{p-1}{2}}-1}{[p]_{q}}.
\endalign
$$
\endproclaim
The proof of Theorem 1.1 will be given in the next sections.

\heading
2. Some Lemmas
\endheading
Below we assume that $p\geq 5$ is a prime and $m$ is an integer prime to $p$.
\proclaim{Lemma 2.1}
$$
\sum_{j=1}^{p-1}\frac{1}{[j]_{q^2}}\equiv-(1+q)\fq(2,q)\pmod{[p]_q}.
$$
\endproclaim
\Proof. This is an immediately consequence of Theorem 1.1 in [P] by observing that
$$
\sum_{j=1}^{p-1}\frac{1}{[j]_{q^2}}=(1+q)\sum_{j=1}^{p-1}\frac{1}{[2j]_{q}}.
$$

\proclaim{Lemma 2.2} Let $m'$ be an integer such that
$$
m'm\equiv 1\pmod{p}.
$$
Then
$$
2\sum_{j\in R_p(m)}\frac{1}{[2j]_{q}}\equiv|R_p(m)|(1-q)+\frac{\fq_p(2,q^{m'})}{[m']_{q}}-\fq_p(2,q)\pmod{[p]_q}.
$$
\endproclaim
\Proof. Clearly
$$
\align
R_p(-m)=&\{1\leq j<p/2:\, \ab{-jm}_p>p/2\}\\
=&\{1\leq j<p/2:\, \ab{jm}_p<p/2\}\\
=&\{1,2,\ldots,(p-1)/2\}\setminus R_p(m).\tag 2.1
\endalign
$$
So applying Lemma 2.1,
$$
\sum_{j\in R_p(m)}\frac{1}{[j]_{q^2}}+\sum_{j\in R_p(-m)}\frac{1}{[j]_{q^2}}=\sum_{j=1}^{(p-1)/2}\frac{1}{[j]_{q^2}}\equiv -(1+q)\fq_p(2,q)\pmod{[p]_q}.\tag 2.2
$$
On the other hand, note that
$\ab{jm}_p>p/2$ if and only if there exists a $1\leq k<p/2$ such that
$$
jm\equiv -k\pmod{p},
$$
or equivalently
$$
j\equiv -km'\pmod{p}.
$$
It follows that
$$
R_p(m)=\{\ab{-km'}_p:\,k\in R_p(m')\},
$$
whence
$$
\sum_{j\in R_p(m)}\frac{1}{[j]_{q^2}}\equiv\sum_{k\in R_p(m')}\frac{1}{[-km']_{q^2}}\pmod{[p]_{q^2}}.
$$
Thus
$$
\align
&\sum_{j\in R_p(m)}\frac{1}{[j]_{q^2}}-\sum_{j\in R_p(-m)}\frac{1}{[j]_{q^2}}\\
\equiv&\sum_{k\in R_p(m')}\frac{1}{[-km']_{q^2}}-\sum_{k\in R_p(-m')}\frac{1}{[km']_{q^2}}\\
=&-\sum_{k\in R_p(m')}\frac{q^{2km'}}{[km']_{q^2}}-\sum_{k\in R_p(-m')}\frac{1}{[km']_{q^2}}\\
=&(1-q^2)|R_p(m')|-\sum_{k\in R_p(m')}\frac{1}{[m']_{q^2}[k]_{q^{2m'}}}-\sum_{k\in R_p(-m')}\frac{1}{[m']_{q^2}[k]_{q^{2m'}}}\\
\equiv&(1-q^2)|R_p(m)|+\frac{(1+q^{m'})\fq_p(2,q^{m'})}{[m']_{q^2}}\pmod{[p]_q}.\tag 2.3
\endalign
$$
Adding (2.2) and (2.3), we obtain that
$$
\sum_{j\in R_p(m)}\frac{1}{[j]_{q^2}}
\equiv\frac{|R_p(m)|}{2}(1-q^2)+\frac{(1+q^{m'})\fq_p(2,q^{m'})}{2[m']_{q^2}}-\frac{1+q}{2}\fq_p(2,q)\pmod{[p]_q}.\tag 2.4
$$
Finally,
$$
\align
2\sum_{j\in R_p(m)}\frac{1}{[2j]_{q}}=&\frac{2}{1+q}\sum_{j\in R_p(m)}\frac{1}{[j]_{q^2}}\\
\equiv&|R_p(m)|(1-q)+\frac{\fq_p(2,q^{m'})}{[m']_{q}}-\fq_p(2,q)\pmod{[p]_q}.
\endalign
$$
\qed
\proclaim{Lemma 2.3}
$$
\align
&\frac{\jacob{m}{p}q^{2\sum_{j\in R_p(m)}(p-\ab{jm}_p)}(q^{2m};q^{2m})_{\frac{p-1}{2}}\big/(q^2;q^2)_{\frac{p-1}{2}}-1}{[p]_{q}}\\
\equiv&2\sum_{j=1}^{(p-1)/2}\left\lfloor\frac{jm}{p}\right\rfloor\frac{q^{2jm}}{[2jm]_{q}}+\fq_p(2,q)-\frac{\fq_p(2,q^m)}{[m]_{q}}-|R_p(m)|(1-q)\pmod{[p]_q}.
\endalign
$$
\endproclaim
\Proof.
Let $r_j=\ab{jm}_{p}$ for any $j\in\Z$. From (2.1), we have
$$
\frac{(q^{2m};q^{2m})_{(p-1)/2}}{(q^2;q^2)_{(p-1)/2}}=\prod_{j=1}^{(p-1)/2}\frac{[jm]_{q^2}}{[j]_{q^2}}=\prod_{j\in R_p(m)}\frac{[jm]_{q^2}}{[j]_{q^2}}\prod_{j\in R_p(-m)}\frac{[jm]_{q^2}}{[j]_{q^2}}.
$$
Now
$$
\prod_{j\in R_p(-m)}\frac{[jm]_{q^2}}{[j]_{q^2}}=\prod_{j\in R_p(-m)}\frac{1-q^{2r_j}}{1-q^{2j}}\(1+\frac{q^{2r_j}(1-q^{2\lfloor jm/p\rfloor p})}{1-q^{2r_j}}\),
$$
and
$$
\prod_{j\in R_p(m)}\frac{[jm]_{q^2}}{[j]_{q^2}}=\prod_{j\in R_p(m)}\frac{1-q^{2(r_j-p)}}{1-q^{2j}}\(1+\frac{q^{2(r_j-p)}(1-q^{2(\lfloor jm/p\rfloor+1)p})}{1-q^{2(r_j-p)}}\).
$$
It is easy to check that
$$
\{\ab{-jm}_p:\,j\in R_p(m)\}\cup\{\ab{jm}_p:\,j\in R_p(-m)\}=\{1,2,\ldots,(p-1)/2\}.
$$
Then
$$
\align
&\prod_{j\in R_p(m)}\frac{1-q^{2(r_j-p)}}{1-q^{2j}}\prod_{j\in R_p(-m)}\frac{1-q^{2r_j}}{1-q^{2j}}\\
=&(-1)^{|R_p(m)|}q^{-2\sum_{j\in R_p(m)}(p-r_j)}\prod_{j\in R_p(m)}\frac{1-q^{2(p-r_j)}}{1-q^{2j}}\prod_{j\in R_p(-m)}\frac{1-q^{2r_j}}{1-q^{2j}}\\
=&\jacob{m}{p}q^{-2\sum_{j\in R_p(m)}(p-r_j)}.
\endalign
$$
Hence
$$
\align
&\jacob{m}{p}q^{2\sum_{j\in R_p(m)}(p-r_j)}\prod_{j=1}^{(p-1)/2}\frac{[jm]_{q^2}}{[j]_{q^2}}\\
\equiv&\prod_{j\in R_p(m)}\(1+\frac{q^{2(r_j-p)}(1-q^{2(\lfloor jm/p\rfloor+1)p})}{1-q^{2(r_j-p)}}\)\prod_{j\in R_p(-m)}\(1+\frac{q^{2r_j}(1-q^{2\lfloor jm/p\rfloor p})}{1-q^{2r_j}}\)\\
\equiv&\prod_{j\in R_p(m)}\(1+\frac{q^{2jm}(1-q^{2(\lfloor jm/p\rfloor+1)p})}{1-q^{2jm}}\)\prod_{j\in R_p(-m)}\(1+\frac{q^{2jm}(1-q^{2\lfloor jm/p\rfloor p})}{1-q^{2jm}}\)\\
\equiv&1+[p]_{q^2}\sum_{j=1}^{(p-1)/2}\left\lfloor\frac{jm}{p}\right\rfloor\frac{q^{2jm}}{[jm]_{q^2}}+[p]_{q^2}\sum_{j\in R_p(m)}\frac{q^{2jm}}{[jm]_{q^2}}\pmod{[p]_{q^2}^2},
\endalign
$$
where in the last setup we use the congruence
$$
\frac{1-q^{jp}}{1-q^p}=1+q^p+\cdots+q^{(j-1)p}\equiv j\pmod{[p]_q}.
$$
Applying (2.4), we have
$$
\align
&\sum_{j\in R_p(m)}\frac{q^{2jm}}{[jm]_{q^2}}\\
=&\sum_{j\in R_p(m)}\frac{1}{[m]_{q^2}[j]_{q^{2m}}}-(1-q^2)|R_p(m)|\\
\equiv&\frac{(1+q^{m'm})\fq_p(2,q^{m'm})}{2[m]_{q^2}[m']_{q^{2m}}}-\frac{1+q^m}{2[m]_{q^2}}\fq_p(2,q^m)+\(\frac{1-q^{2m}}{2[m]_{q^2}}-(1-q^2)\)|R_p(m)|\\
\equiv&\frac{1+q}{2}\fq_p(2,q)-\frac{1+q^m}{2[m]_{q^2}}\fq_p(2,q^m)-\frac{|R_p(m)|}{2}(1-q^2)\pmod{[p]_q}.
\endalign
$$
Thus
$$
\align
&\frac{\jacob{m}{p}q^{2\sum_{j\in R_p(m)}(p-r_j)}(q^{2m};q^{2m})_{\frac{p-1}{2}}\big/(q^2;q^2)_{\frac{p-1}{2}}-1}{[p]_{q}}\\
=&\frac{1+q^p}{1+q}\cdot\frac{\jacob{m}{p}q^{2\sum_{j\in R_p(m)}(p-r_j)}(q^{2m};q^{2m})_{\frac{p-1}{2}}\big/(q^2;q^2)_{\frac{p-1}{2}}-1}{[p]_{q^2}}\\
\equiv&\frac{2}{1+q}\sum_{j=1}^{(p-1)/2}\left\lfloor\frac{jm}{p}\right\rfloor\frac{q^{2jm}}{[jm]_{q^2}}+\frac{2}{1+q}\sum_{j\in R_p(m)}\frac{q^{2jm}}{[jm]_{q^2}}\\
\equiv&2\sum_{j=1}^{(p-1)/2}\left\lfloor\frac{jm}{p}\right\rfloor\frac{q^{2jm}}{[2jm]_{q}}+\fq_p(2,q)-\frac{\fq_p(2,q^m)}{[m]_{q}}-|R_p(m)|(1-q)\pmod{[p]_q}.
\endalign
$$
\qed
\heading
3. Proof of Theorem 1.1
\endheading
We write
$$
\align
\qbinom{p-1}{\floor{kp/m}}{q^{2m}}=&\prod_{j=1}^{\floor{kp/m}}\frac{[p-j]_{q^{2m}}}{[j]_{q^{2m}}}\\
=&\prod_{j=1}^{\floor{kp/m}}\frac{[p]_{q^{2m}}-[j]_{q^{2m}}}{q^{2jm}[j]_{q^{2m}}}\\
=&(-1)^{\floor{kp/m}}q^{-2m\binom{\floor{kp/m}+1}{2}}\prod_{j=1}^{\floor{kp/m}}\(1-\frac{[p]_{q^{2m}}}{[j]_{q^{2m}}}\).
\endalign
$$
Then
$$
\align
&(-1)^{\sum_{k=1}^{\floor{m/2}}\floor{kp/m}}q^{2m\sum_{k=1}^{\floor{m/2}}\binom{\floor{kp/m}+1}{2}}\prod_{k=1}^{(m-1)/2}\qbinom{p-1}{\floor{kp/m}}{q^{2m}}\\
\equiv&1-[p]_{q^{2m}}\sum_{k=1}^{\floor{m/2}}\sum_{j=1}^{\floor{kp/m}}\frac{1}{[j]_{q^{2m}}}\\
\equiv&1-[p]_{q^{2m}}\sum_{j=1}^{(p-1)/2}\sum_{k=\ceil{jm/p}}^{\floor{m/2}}\frac{1}{[j]_{q^{2m}}}\\
\equiv&1-[p]_{q^{2m}}\sum_{j=1}^{(p-1)/2}\frac{\floor{m/2}-\floor{jm/p}}{[j]_{q^{2m}}}\pmod{[p]_{q^{2m}}^2}.
\endalign
$$
By Lemma 2.1,
$$
\align
[p]_{q^{2m}}\sum_{j=1}^{(p-1)/2}\frac{\floor{m/2}}{[j]_{q^{2m}}}\equiv&-\floor{\frac{m}{2}}[p]_{q^{2m}}(1+q^m)\fq_p(2,q^m)\\
=&-\floor{\frac{m}{2}}\frac{1-q^{2mp}}{1-q^m}\fq_p(2,q^m)\\
\equiv&-2\floor{m/2}[p]_{q^m}\fq_p(2,q^m)\pmod{[p]_{q^m}^2}.
\endalign
$$
And
$$
\align
[p]_{q^{2m}}\sum_{j=1}^{(p-1)/2}\frac{\floor{jm/p}}{[j]_{q^{2m}}}=&[mp]_{q^{2}}\sum_{j=1}^{(p-1)/2}\frac{\floor{jm/p}}{[jm]_{q^{2}}}\\
=&[mp]_{q^{2}}\sum_{j=1}^{(p-1)/2}\floor{\frac{jm}{p}}\frac{q^{2jm}}{[jm]_{q^{2}}}+(1-q^{2mp})\sum_{j=1}^{(p-1)/2}\floor{\frac{jm}{p}}.
\endalign
$$
From Lemma 2.3, we deduce that
$$
\align
&[mp]_{q^{2}}\sum_{j=1}^{(p-1)/2}\floor{\frac{jm}{p}}\frac{q^{2jm}}{[jm]_{q^{2}}}\\
\equiv&2m[p]_q\sum_{j=1}^{(p-1)/2}\floor{\frac{jm}{p}}\frac{q^{2jm}}{[2jm]_{q}}\\
\equiv&m[p]_q\(\eq_p^*(m,q)-\fq_p(2,q)+\frac{\fq_p(2,q^m)}{[m]_{q}}+|R_p(m)|(1-q)\)\pmod{[p]_q^2}.
\endalign
$$
Since $p\nmid 2m$, we have both $[p]_{q^m}$ and $[p]_{q^{2m}}$ are divisible by $[p]_q$.
Also note that
$$
[p]_{q^m}=[mp]_{q}/[m]_q\equiv m[p]_q/[m]_q\pmod{[p]_q^2}
$$
and
$$
1-q^{2mp}\equiv 2m(1-q)[p]_q\pmod{[p]_q^2}.
$$
Therefore we obtain that
$$
\align
&(-1)^{\sum_{k=1}^{(m-1)/2}\floor{\frac{kp}{m}}}q^{2m\sum_{k=1}^{(m-1)/2}\binom{\floor{kp/m}+1}{2}}\prod_{k=1}^{(m-1)/2}\qbinom{p-1}{\floor{kp/m}}{q^{2m}}\\
\equiv&1+m[p]_q\eq_p^*(m,q)+(2\floor{m/2}+1)[p]_{q^m}\fq_p(2,q^m)-m[p]_q\fq_p(2,q)\\
&+m\(|R_p(m)|+2\sum_{j=1}^{(p-1)/2}\floor{\frac{jm}{p}}\)(1-q)[p]_q\pmod{[p]_q^2}.
\endalign
$$
Finally, by Lemma 3.1 in [S],
$$
\sum_{k=1}^{(m-1)/2}\floor{\frac{kp}{m}}\equiv \frac{p-1}{2}\floor{\frac{m}{2}}+\frac{(p^2-1)(m-1)}{8}-|R_p(m)|\pmod{2},
$$
which implies that
$$
(-1)^{\sum_{k=1}^{(m-1)/2}\floor{\frac{kp}{m}}}=(-1)^{\frac{p-1}{2}\floor{\frac{m}{2}}}\jacob{m}{p}\jacob{2}{p}^{m-1}.
$$
All are done.\qed

\Ack. We thank our advisor, Professor Zhi-Wei Sun, for his help on this paper.

\enddocument